\documentclass[11pt]{article}
\usepackage{amsmath,amsthm,amsfonts,comment}
\DeclareMathOperator{\1}{id}

\newcommand{\p}{\partial}
\renewcommand{\=}{\doteq}

\newtheorem{thm}{Theorem}
\newtheorem{cor}[thm]{Corollary}
\theoremstyle{definition}

\theoremstyle{definition}
 
\allowdisplaybreaks

\begin{document}

\thispagestyle{empty} 

\title{\bf\Large Moufang loops and generalized Lie-Cartan theorem\footnote{Presented at the International Conference  "Loops '07", Prague, August 12-19, 2007.}} 
\author{\Large Eugen Paal\\ \\
Department of Mathematics\\
Tallinn University of Technology\\
Ehitajate tee 5, 19086 Tallinn, Estonia\\ 
\smallskip
E-mail: eugen.paal@ttu.ee}
\date{}

\maketitle

\begin{abstract}
Generalized Lie-Cartan theorem for linear birepresentations of an analytic Moufang loop is considered.
The commutation relations of the generators of the birepresentation are found. 
In particular, the Lie algebra of the multiplication group of the birepresentation is explicitly given.
\end{abstract}

\begin{flushright}
\emph{Dedicated to Maks A. Akivis on the occasion of his\\ 85th birthday and 65 years of scientific activity}
\end{flushright}
\section{Introduction}

Continuous symmetries generated by the Lie transformation groups are widely exploited in modern mathematics and physics. Nevertheless, it may  happen that the group theoretical methods are too rigid and one has to extend these beyond the Lie groups and algebras. From this point of view it is interesting to extend the group theoretical symmetry methods by using the Moufang loops and Mal'tsev algebras. The latter are known as minimal non-associative generalizations of the group and Lie algebra concepts, respectively. 

In this paper, the generalized Lie-Cartan theorem for linear birepresentations of an analytic Moufang loop is considered.  The commutation relations of the generators of the birepresentation were found. In particular, the Lie algebra of the multiplication group of the birepresentation is explicitly given.

Based on this theorem, various applications are possible. 
In particular, it was recently shown \cite{Paal06} how the Moufang-Noether current algebras may be constructed so that the corresponding Noether charge algebra turns out to be a birepresentation of the tangent Mal'tsev algebra of an analytic Moufang loop. 

\section{Moufang loops}

A {\it Moufang loop} \cite{RM,Bel,HP} is a quasigroup $G$ with the unit
element $e\in G$ and the Moufang identity
\[
(ag)(ha)=a(gh)a,\qquad a,g,h\in G
\]
Here the multiplication is denoted by juxtaposition.
In general, the multiplication need not be associative: $gh\cdot a\neq g\cdot ha$ for some triple of elements 
$a,g,h\in G$. The inverse element $g^{-1}$ of $g$ is defined by
\[
gg^{-1}=g^{-1}g=e
\]

\section{Analytic Moufang loops and Mal'tsev algebras}

Following the concept of the Lie group, the notion of an  analytic
Moufang loop can be easily formulated.

A Moufang loop $G$ is said \cite{Mal} to be {\it analytic} if $G$ is also a real
analytic manifold and main operations - multiplication and inversion map
$g\mapsto g^{-1}$ - are analytic mappings.

Let the local coordinates of $g$ from the vicinity of $e\in G$ be
denoted by $g^{i}$ ($i=1,\dots,r\=\dim G$). The \emph{tangent space} of
$G$ at $e\in G$ is denoted by $T_e(G)$.

As in the case of the Lie groups, the structure constants $c^{i}_{jk}$
of an analytic Moufang loop are defined by
\[
c^{i}_{jk}\=\frac{\p^{2}(ghg^{-1}h^{-1})^{i}}{\p g^{j}\p h^{k}}\Big|_{g=h=e}=-c^{i}_{kj},
\qquad i,j,k=1,\dots,r
\]
For any $x,y\in T_e(G)$, their (tangent)
product $[x,y]\in T_e(G)$ is
defined in component form by
\[
[x,y]^{i}\=c^{i}_{jk}x^{j}y^{k}=-[y,x]^{i},\qquad i=1,\dots,r
\]
The tangent space $T_e(G)$ being equipped with such an anti-commutative
multiplication is called the {\it tangent algebra} of the analytic Moufang loop $G$.
We shall use the notation $\Gamma\=\{T_e(G),[\cdot,\cdot]\}$ for the tangent algebra of $G$.

The tangent algebra of $G$ need not be a Lie algebra. There may exist a triple
$x,y,z\in T_e(G)$ that does not satisfy the Jacobi identity:
\[
J(x,y,z)\=[x,[y,z]]+[y,[z,x]]+[z,[x,y]]\neq0
\]
Instead, for all $x,y,z\in T_e(G)$ one has a more general {\it Mal'tsev identity}
\cite{Mal}
\[
[J(x,y,z),x]=J(x,y,[x,z])
\]
Anti-commutative algebras with this identity are called the {\it Mal'tsev algebras}.

\section{Birepresentation of a Moufang loop}

Consider a pair $(S,T)$ of the maps $g\mapsto S_g$, $g\mapsto T_g$  of a Moufang loop $G$ into $GL_n$. The pair $(S,T)$ is called \cite{Paal98} a (linear)  \emph{birepresentation} of $G$ (in $GL_n)$  if the following conditions hold true:

\begin{itemize}
\itemsep-1pt
\item
$S_e=T_e=\1$,
\item
$T_gS_gS_h=S_{gh}T_g$,
\item
$S_gT_gT_h=T_{hg}S_g$.
\end{itemize}

The birepresentation $(S,T)$ is called \emph{associative}, if
the following simultaneous relations are satisfied:
\[
S_gS_h=S_{gh},\quad T_gT_h=T_{hg},\quad S_gT_h=T_hS_g,\qquad \forall\, g,h\in G
\]
In general, birepresentations need not be associative even for groups.

\section{ Generators and structure functions}

The \emph{generators} of a birepresentation $(S,T)$ are defined as follows:
\[
S_j\=\frac{\p S_g}{\p g^{j}}\Bigg\vert_{g=e},
\qquad
T_j\=\frac{\p T_g}{\p g^{j}}\Bigg\vert_{g=e},
\qquad
j=1,\dots,r
\]
The \emph{auxiliary functions} of $G$ are defined as
\[
v_j^{n}(g)\=\frac{\p(gh)^{n}}{\p h^{j}}\Bigg\vert_{h=e},
\qquad
n,j=1,\dots,r
\]
The \emph{structure functions} $c_{jk}^{n}(g)$ of $G$
are defined by the \emph{generalized Maurer-Cartan equations}
\[
v_j^{n}(g)\frac{\p v_k^{i}(g)}{\p g^{n}}-v_k^{n}(g)\frac{\p v_j^{i}(g)}{\p g^{n}}
\=c_{jk}^{n}(g)v_n^{i}(g)
\]
One can check the initial conditions
\[
c_{jk}^{n}(e)=c_{jk}^{n}
\]
The algebra with structure functions $c_{jk}^{n}(g)$ may be called the \emph{derivative} of the tangent algebra $\Gamma$ of $G$. It follows from the Belousov theory \cite{Bel} of \emph{derivative quasigroups and loops} that  the derivative of the tangent Mal'tsev algebra of $G$  is a Mal'tsev algebra as well. In a sense, the derivative Mal'tsev algebra is a deformation of $\Gamma$ with the deformation parameter $g\in G$.

\section{Generalized Lie equations (GLE)}

Define the \emph{derivative generators} of $(S,T)$ as follows:
\[
S_j(g)\=T_gS_jT_g^{-1}, \quad T_j(g)\=S_g^{-1}T_jS_g
\]
By direct calculations, one can check that
\[
S_j(g)+T_j(g)=S_j+T_j
\]
In what follows, the \emph{commutator} of linear operators $A,B$ is denoted as $[A,B]\=AB-BA$.

\begin{thm}[generalized Lie equations]
Let $(S,T)$ be a differentiable birepresentation of an analytic
Moufang loop $G$. Then the birepresentation matrices $S_g,T_g$
($g\in G$) satisfy the system of simultaneus differential
equations 
\begin{align*}
v_j^{n}(g)\frac{\p S_g}{\p g^{n}}
&=S_g\,\overbrace{T_gS_jT_g^{-1}}^{S_j(g)}=S_g\,S_j(g)\\
&=S_gS_j+\underbrace{[S_g,T_j]}_{\text{associator}}\\
v_j^{n}(g)\frac{\p T_g}{\p g^{n}}
&=\overbrace{S_g^{-1}T_jS_g}^{T_j(g)}\,T_g=T_j(g)\,T_g\\
&=T_jT_g+\underbrace{[S_j,T_g]}_{\text{associator}}.
\end{align*}
\end{thm}
\begin{proof}
Differentiate the defining relations of a birepresentation $(S,T)$ with
respect to $g$  and take $g=e$.
\end{proof}
\section{Generalized Lie-Cartan theorem}

\begin{thm}[generalized Lie-Cartan theorem]
The integrability conditions of GLE read as commutation relations
\begin{align*}
[S_j(g),S_k(g)]&=\hphantom{-}c_{jk}^{n}(g)S_n(g)-2[S_j(g),T_k(g)]\\
[T_j(g),T_k(g)]&=-c_{jk}^{n}(g)T_n(g)-2[T_j(g),S_k(g)]
\end{align*}
\end{thm}
\begin{proof}
Differentiate the GLE and use
\[
\frac{\p^{2}S_g}{\p g^{j}\p g^{k}}=\frac{\p^{2}S_g}{\p g^{k}\p g^{j}}\,,
\qquad
\frac{\p^{2}T_g}{\p g^{j}\p g^{k}}=\frac{\p^{2}T_g}{\p g^{k}\p g^{j}}
\tag*{\qed}
\]
\renewcommand{\qed}{}
\end{proof}

\begin{cor}
If $(S,T)$ is an associative birepresentation of a Lie group $G$, then one gets the well known Lie algebra
commutation relations (Lie-Cartan theorem)
\[
[S_j,S_k]-c_{jk}^{n}S_n
=[T_j,T_k]+c_{jk}^{n}T_n
=[S_j,T_k]
=0
\]
\end{cor}
\section{Yamagutian and Yamaguti brackets}

Let $g\in G$ and $x,y\in T_e(G)$. Denote
\begin{gather*}
S_x(g)\=x^{j}S_j(g),\qquad T_x(g)\=x^{j}T_j(g),\qquad
[x,y]^{i}_g\=c^{i}_{jk}(g)x^{j}y^{k},\quad i=1,\dots,r
\end{gather*}
Define the \emph{Yamagutian} $Y_g$ by
\[
3Y_g(x;y)\=[S_x(g),S_y(g)]+[S_x(g),T_y(g)]+[T_x(g),T_y(g)]
\]
By direct calulations one can check that the Yamagutian obeys the constraints
\begin{gather*}
Y_g(x;y)+Y_g(y;x)=0\\
Y_g([x,y]_g;z)+Y_g([y,z]_g;x)+Y_g([z,x]_g;y)=0
\end{gather*}
The \emph{Yamaguti brackets} $[\cdot,\cdot,\cdot]_g$ are defined \cite{Yam62, Yam63} in $T_e(G)$ by
\begin{align*}
[x,y,z]_g
&\=[x,[y,z]_g]_g-[y,[x,z]_g]_g+[[x,y]_g,z]_g\\
&=J_g(x,y,z)+2[[x,y]_g,z]_g
\end{align*}

\section{Closure of integrability conditions}

It turns out that the commutation relations for derivative generators of the
birepresentation $(S,T)$ can be presented in a closed Lie algebra form. 

The integrability conditions of GLE can be written as follows:
\begin{align*}
[S_x(g),S_y(g)]&=\,\,2Y_g(x;y)+\frac{1}{3}S_{[x,y]_g}(g)+\frac{2}{3}T_{[x,y]_g}(g) \tag{1}\\
[S_x(g),T_y(g)]&=-Y_g(x;y)+\frac{1}{3}S_{[x,y]_g}(g)-\frac{1}{3}T_{[x,y]_g}(g) \tag{2}\\
[T_x(g),T_y(g)]&=\,\,2Y_g(x;y)-\frac{2}{3}S_{[x,y]_g}(g)-\frac{1}{3}T_{[x,y]_g}(g) \tag{3}
\end{align*}
One can prove  the \emph{reductivity conditions} \cite{Paal98}
\begin{align*}
6[Y_g(x;y),S_z(g)]&=S_{[x,y,z]_g}(g) \tag{4}\\
6[Y_g(x;y),T_z(g)]&=T_{[x,y,z]_g}(g) \tag{5}
\end{align*}
By using the above reductivity conditions it can be shown by direct calculations 
that the \hbox{Yamagutian} obeys the Lie algebra
\begin{align*}
6[Y_g(x;y),Y_g(z;w)]&=Y_g([x,y,z]_g;w)+Y_g(z;[x,y,w]_g) \tag{6}
\end{align*}
Dimension of the Lie algebra $(1)-(6)$ does not exceed $2r+r(r-1)/2$. The
Jacobi identities are guaranteed by the defining identities of the \emph{Lie} \cite{Loos}
and \emph{general Lie triple systems} \cite{Yam62,Yam63} associated with the derivative Mal'tsev
algebra of $T_e(G)$ of $G$.

The commutation relations of form $(1)-(6)$ are well known from the theory of \emph{alternative} algebras \cite{Schafer}. By taking $g=e$, the commutation relations $(1)-(6)$ give  the Lie algebra of the 
\emph{multiplication group of the birepresentation} $(S,T)$ of $G$. For field theoretical applications see  \cite{Paal06}.

\section*{Acknowledgment}
The research was in part supported by the Estonian  Science Foundation, Grant 6912.

\end{document}